\documentclass[10pt]{article}
\usepackage[utf8]{inputenc}
\usepackage{epsfig}
\usepackage{color}
\usepackage[british,english]{babel}
\usepackage{amsthm}
\usepackage{amsmath}
\usepackage{amsfonts}
\usepackage{amssymb}
\usepackage{graphicx}
\newtheorem{corollary}{Corollary}[section]

\newtheorem{lemma}[corollary]{Lemma}

\newtheorem{remark}[corollary]{Remark}
\newtheorem{theorem}[corollary]{Theorem}
\newfont{\sBlackboard}{msbm10 scaled 900}

\newcommand{\mylabel}[1]{\label{#1}
            \ifx\undefined\stillediting
            \else \fbox{$#1$}\fi }
\newcommand{\BE}{\begin{equation}}

\newcommand{\EEQ}{\end{equation}}
\newcommand{\rfb}[1]{\mbox{\rm
   (\ref{#1})}\ifx\undefined\stillediting\else:\fbox{$#1$}\fi}

\newfont{\Blackboard}{msbm10 scaled 1200}

\newfont{\roma}{cmr10 scaled 1200}

\def\CC{\rm \hbox{C\kern-.56em\raise.4ex
         \hbox{$\scriptscriptstyle |$}\kern+0.5 em }}

\def\n{|\kern -.05cm{|}\kern -.05cm{|}}

\def \noame{\noalign{\medskip}}
\newcommand{\mm}    {{\hbox{\hskip 0.5pt}}}

\newcommand{\bluff} {{\hbox{\raise 15pt \hbox{\mm}}}}

\newcommand{\ep}   {\varepsilon}
\usepackage{fancyhdr}

\makeatletter
\def\section{\@startsection {section}{1}{\z@}{-3.5ex plus -1ex minus
    -.2ex}{2.3ex plus .2ex}{\large\bf}}
\makeatother
\def\be{\begin{equation}}
\def\ee{\end{equation}}

\date{ }
\begin{document}
\thispagestyle{empty}
\title{\Large \bf  Mathematical modelling of a thin-film flow  obeying  Carreau's law  without high-rate viscosity}\maketitle
\vspace{-2cm}
\begin{center}
Mar\'ia ANGUIANO\footnote{Departamento de An\'alisis Matem\'atico. Facultad de Matem\'aticas. Universidad de Sevilla. 41012-Sevilla (Spain) anguiano@us.es} and Francisco Javier SU\'AREZ-GRAU\footnote{Departamento de Ecuaciones Diferenciales y An\'alisis Num\'erico. Facultad de Matem\'aticas. Universidad de Sevilla. 41012-Sevilla (Spain) fjsgrau@us.es}
 \end{center}

 \renewcommand{\abstractname} {\bf Abstract}
\begin{abstract}
 In this paper, we derive an extension of the Reynolds law for  quasi-Newtonian fluid flows through a thin domain with thickness $0<\varepsilon\ll 1$ with viscosity obeying the Carreau law without high-rate viscosity, by applying asymptotic analysis with respect to $\varepsilon$. This provides a framework for understanding how the  non-Newtonian effects and the thickness of the domain (which is significantly smaller than the other dimensions) influence its flow behavior. 
\end{abstract}
\bigskip\noindent

\noindent {\small \bf AMS classification numbers:}  35B27; 35Q35; 76A05; 76M50; 76A20.  \\

\ \\
\section {Introduction}\label{S1}
The Reynolds equation is a fundamental equation in lubrication theory that describes the pressure distribution in a thin film of fluid between two surface. It is often used in the study of bearings, seals, and other machine elements where a thin fluid layer separates moving parts.  Using the film thickness as a small parameter, an asymptotic approximation of the Stokes system can be derived providing the well-known Reynolds equation for the pressure $\widetilde p$ of the fluid, see Bayada \& Chambat \cite{Bayada_Chambat_Transition} for more details
 \begin{equation}\label{ReynoldsNew}
 {\rm div}_{y'}\left({h(y')^3\over 12\mu}\left({ f'}(y')-\nabla_{y'}\widetilde p(y')\right)\right)=0, 
 \end{equation}
 where $h$ describes the shape of the top boundary,  $\mu$ is the viscosity of the fluid, $f'=(f_1,f_2)$ is the body forces and $y'\in\omega\subset\mathbb{R}^2$. We refer to Bayda \& Chambat \cite{Bayada2} for the homogenization of the Stokes system in a thin film flow with rapidly varying thickness, to Bresch {\it et al.} \cite{Bresch} for the study of the roughness-induced effect  on the Reynolds equation,  and  Duvnjak \& Maru${\rm \check{s}}$i\'c-Paloka \cite{Duvnjak2} for the derivation of the Reynolds equation for lubrication of a rotating shaft. 

However, many others fluids behave differently, the viscosity of these fluids is no more constant, such as pastes or polymer solutions. These fluids  are called non-Newtonian. The simplest idea to describe non-Newtonian fluids is to plot the viscosity measurements versus the imposed shear rate and then, to fit the obtained curve with a simple template viscosity function, adjusting some few parameters. This is the main idea of quasi-Newtonian fluids models, which could be viewed as a first step inside the world of non-Newtonian fluids models (see Saramito \cite[Chapter 2]{Saramito} for more details).

 The incompressible quasi-Newtonian fluids are characterized by the viscosity depending on
 the principal invariants of the symmetric stretching tensor $\mathbb{D}[{ u}]$. If ${\  u}$ is the velocity, $p$ is the pressure and $D{ u}$ is the gradient velocity tensor, $\mathbb{D}[{ u}]=(D{ u}+D^t {  u})/2$ denotes the symmetric stretching tensor and $\sigma$ the stress tensor given by $\sigma=-pI+2\eta_r \mathbb{D}[{\  u}]$. The viscosity $\eta_r$ is constant for a Newtonian fluid but dependent on the shear rate, that is, $\eta_r=\eta_r(\mathbb{D}[{  u}])$, for viscous non-Newtonian fluids. The deviatoric stress tensor $\tau$, that is, the part of the total stress tensor that is zero at equilibrium, is then a nonlinear function of the shear rate $\mathbb{D}[{  u}]$, that is $\tau=\eta_r(\mathbb{D}[{\  u}])\mathbb{D}[{  u}]$, see Baranger \& Najib \cite{Baranger}, Barnes {\it et al.} \cite{Barnes}, Bird {\it et al.} \cite{Bird},  Mikeli\'c \cite{MikelicIntro2} and Sandri \cite{Sandri} for more details).
 
 Two most widely used laws in engineering practice are the power law and the Carreau law.   The most popular model is the power law where the expression for the shear rate dependent viscosity is 
$$\eta_r(\mathbb{D}[{\bf u}])=\mu |\mathbb{D}[{\bf u}]|^{r-2},\quad 1<r<+\infty,$$
where the two material parameters $\mu>0$ and $r$ are called the consistency and the flow index, respectively. Here, the matrix norm $|\cdot |$ is defined by $|\xi|^2=Tr(\xi\xi^t)$ with $\xi\in \mathbb{R}^3$. We recall that   $1<r<2$ corresponds to a pseudoplastic fluid (shear thinning), while for $r>2$ the fluid is dilatant (shear thickening). We remark that $r=2$ yields the Newtonian fluid. Similarly to the derivation of the classical Reynolds equation, a two-dimensional nonlinear Reynolds equation for quasi-Newtonian fluid with viscosity given by the power law has been obtained in Mikeli\'c \& Tapiero \cite{Tapiero}, with the following form
\begin{equation}\label{ReynoldsPowerLaw}
{\rm div}_{y'}\left(
{h(y')^{r'+1}\over 2^{r'\over 2}(r'+1)\mu^{r'-1}}\left|{ f}'(y')-\nabla_{y'}\widetilde p(y')\right|^{r'-2}\left({ f}'(y')-\nabla_{y'}\widetilde p(y')\right)\right)=0,\quad\hbox{in }\omega,
\end{equation}
where $r'$ is the conjugate exponent of $r$, i.e. $1/r+1/r'=1$. Also, we refer to Anguiano \& Su\'arez-Grau \cite{Anguiano_nonlinear} for the derivation of nonlinear Reynolds equations for non-Newtonian thin-film fluid flows over a rough boundary, to Duvnjak \cite{Duvnjak} for the derivation of non-linear Reynolds-type problem for lubrication of a rotating shaft, and to Fabricius {\it et al.} \cite{Fabricius} for a study on pressure-driven Hele-Shaw flow of power-law fluids.

The simple power-law model has a well-known singularity at zero shear rate, which must be carefully accounted for in kinematic analyses. The Carreau  equation is an alternate generalized Newtonian model that enables the description of the plateaus in viscosity that are expected when the shear rate is very small or very large. The   Carreau  law is given by
\begin{equation}\label{Carreaulaw_entera}\eta_r(\xi)=(\eta_0-\eta_\infty)(1+\lambda|\xi|^2)^{{r\over 2}-1}+\eta_\infty,\quad 1<r<+\infty,\ r\neq 2,\quad\eta_0>\eta_\infty\geq 0,\quad \lambda>0,
\end{equation}
where $\eta_0$ is the zero-shear-rate viscosity, $\eta_\infty$ is the infinite-shear-rate viscosity, $\lambda$ is a time constante being the inverse of a characteristic shear rate at which shear thinning becomes important, $r-1$ is a dimensionless constant describing the slope in the {\it power law region} of $\log(\eta_r)$ versus $\log(|\mathbb{D}[u_\ep]|)$. The case $\eta_\infty>0$ corresponds to a polymer in solution with the solvent's viscosity equal to $\eta_\infty$, and $\eta_\infty=0$ corresponds to a polymer melt.

 For $\eta_\infty>0$,  a two-dimensional Reynolds equation for a  thin film flow of quasi-Newtonian fluids with viscosity obeying the Carreau law was derived  in Boughanim \& Tapiero \cite{Tapiero2} (case $\gamma=1$ and $r\neq 2$), with the following form
\begin{equation}\label{ReynoldsCarreauLaw}
 {\rm div}_{y'}\left(\left(\int_{-{h(y')\over 2}}^{h(y')\over 2}{({h(y')\over 2}+\xi)\xi \over \psi(2|f'(y')-\nabla_{y'}\widetilde p(y')||\xi|)}\,d\xi\right)\left(f'(y')-\nabla_{y'}\widetilde p(y')\right)\right)=0\quad \hbox{in }\omega
\end{equation}
where $\zeta=\psi(\tau)$ for $\tau\in\mathbb{R}^{+}$  is the unique solution noted  of the algebraic equation
\begin{equation}\label{algebraic}
\tau=\zeta\sqrt{{2\over \lambda}\left\{{\zeta-\eta_\infty\over \eta_0-\eta_\infty}\right\}^{2\over r-2}-1}.
\end{equation}
Also, we refer to Bourgeat {\it et al} \cite{Bourgeat_Marusic_Marusic} fo the injection of a non-Newtonian fluid through a thin periodically perforated wall and 
 to Anguiano {\it et al.} \cite{Anguiano_Carreau1, Anguiano_Carreau2} for the modeling of Carreau fluids in thin porous media.

However, to the authors' knowledge, in the mathematical literature there does not appear to be  any derivation of a Reynolds equation through a thin domain for a quasi-Newtonian fluid flow with the Carreau viscosity law with $\eta_\infty=0$ (for a polymer melt), i.e. 
\begin{equation}\label{Carreaulaw}\eta_r(\xi)=\eta_0(1+\lambda|\xi|^2)^{{r\over 2}-1},\quad 1<r<+\infty,\ r\neq 2,\quad\eta_0>0,\quad \lambda>0,
\end{equation}  
which is called the law of Bird and Carreau. This simplification of neglecting the  high-shear viscosity $\eta_\infty$  can also be valid when focusing on the behavior of the fluid at lower shear rates and the transition to a power-law region.  Thus, this simplified Carreau model can be useful in situations where the fluid's behavior at very high shear rates is not of primary interest, or when simpler calculations and analysis are desired while still capturing the essential non-Newtonian characteristics of the fluid.  It is also valid when $\mu_\infty$ is very small, so in the applications it is often set to zero. Considering the modified Carreau law  (\ref{Carreaulaw}), we refer to Bourgeat \& Marusi\'c-Paloka \cite{BourgeatPaloka0} for the mathematical modelling of a quasi-Newtonian viscous flow through a thin filter, and to Mikeli\'c \cite{Mikelicmu0, MikelicIntro2} for the derivation of a Darcy laws  for a quasi-Newtonian fluid flows in a perforated domain.  

In this paper, by using asymptotic analysis,  starting from the three-dimensional Stokes equation with the Carreau law without high rate viscosity $\mu_\infty$,  we derive a two-dimensional Reynolds problem, which agrees with (\ref{ReynoldsCarreauLaw}) where $\psi$ is the  unique solution  of the algebraic equation (\ref{algebraic}) with $\eta_\infty=0$. To prove this result, we first proceed with formal asymptotic expansion of velocity and pressure (see Section  \ref{sec:formal}) and then, by compactness results, we make the rigorous proof (see Theorem \ref{them_limit} and Corollary \ref{cor_thm}). The main mathematical difference with respect to the case $\eta_\infty>0$ is that, in the case $\eta_\infty>0$, the solution (velocity/pressure) of the Stokes system is in $W^{1,r}_0\times L^{r'}_0$ if $r>2$ and in  $H^1_0\times L^{2}_0$ if $1<r<2$, whereas in the present case $\eta_\infty=0$, the solution is in  $W^{1,r}_0\times L^{r'}_0$ if $1<r<+\infty$, $r\neq 2$, so we need to derive the {\it a priori} estimates in a different way with respect to the case $\eta_\infty>0$. In the analysis of the pressure, we introduce a more accurate {\it a priori} estimate by using a decomposition of the pressure in two pressures, one of which is in $W^{1,r'}$ and gives the macroscopic behavior,  that allows to prove strong convergence of the pressure and better regularity of the limit pressure before passing to the limit.  

 Finally, we comment the structure of the paper. In Section \ref{sec:definition} we define the domain and introduce some notation which will be useful in the rest of the paper. In Section \ref{sec:model}, we present the model that will be  analyzed. A formal asymptotic analysis is developed in Section \ref{sec:formal}. A priori estimates of velocity and pressure are derived in Section \ref{sec:estimates}. Convergences of velocity and pressure and the limit model are presented in Section \ref{sec:convergences}. We finish the paper with a list of references.

\section{Definition of the domain and some notation}\label{sec:definition}   We define the thin domain    by
\begin{equation}\label{Omegaep}
\Omega_\varepsilon=\{x=(x',x_3)\in\mathbb{R}^2\times \mathbb{R}\,:\, x'\in \omega,\ 0<x_3<  h_\ep(x')\},
\end{equation}
where  $\omega$ a smooth and connected subset of $\mathbb{R}^2$ and the function $h_\ep(x')= \varepsilon h\left(x'\right)$ represents the distance between the two surfaces, i.e. the thickness of the domain, where $0<\varepsilon\ll 1$. 
\noindent Function $h$ is a positive and $C^1$ function defined for $x'$ in $\omega$, and there exist $h_{\rm min}$ and $h_{\rm max}$ such that
$$0<h_{\rm min} \leq h(x')\leq  h_{\rm max}<+\infty\,.$$

\noindent
We define the boundaries of $\Omega_\varepsilon$ as follows
$$\begin{array}{c}
\displaystyle \Gamma_0 =\omega\times\{0\},\quad \Gamma_1^\varepsilon=\left\{(x',x_3)\in\mathbb{R}^2\times\mathbb{R}\,:\, x'\in \omega,\ x_3=  h_\varepsilon(x')\right\},\quad
\displaystyle\Gamma_{\rm lat}^\varepsilon=\partial\Omega_\ep\setminus (\Gamma_0\cup\Gamma_1^\ep).
\end{array}$$

\noindent Applying a dilatation in the vertical variable, i.e. $y'=x', y_3=x_3/\ep$, we define the following rescaled sets 
 \begin{equation}\label{domains_tilde}\begin{array}{c}
 \displaystyle   \Omega=\{(y',y_3)\in\mathbb{R}^2\times \mathbb{R}\,:\, y'\in \omega,\ 0<y_3< h(y')\},\\
 \noame
 \Gamma_1 =\{(y',y_3)\in\mathbb{R}^2\times \mathbb{R}\,:\, y'\in \omega,\ y_3=h(y')\}\quad   \Gamma_{\rm lat}=\partial  \Omega \setminus (\Gamma_0\cup  \Gamma_1).
  \end{array}
  \end{equation}

\noindent In the sequel, we introduce the following notation. Let us consider a vectorial function ${  v} =( {  v} ', v_{3})$ with ${ v}'=(v_1, v_2)$ and $\varphi$ a scalar function, both defined in $\Omega_\varepsilon$. Then, introduce the operator $\Delta$ and $\nabla$ by
$$\begin{array}{c}
 \displaystyle \Delta{  v}=\Delta_{x'}{  v} +\partial_{x_3}^2 { v},\quad {\rm div}({  v} )={\rm div}_{x'}({  v}')+\partial_{x_3}v_3,\quad \nabla \varphi=(\nabla_{x'} \varphi,  \partial_{x_3}\varphi)^t,\\
 \end{array}$$
and
$\mathbb{D}:\mathbb{R}^3\to \mathbb{R}^3_{\rm sym}$   the symmetric part of the velocity gradient, that is
$$\mathbb{D}[{  v}]={1\over 2}(D{  v}+(D{  v})^t)=\left(\begin{array}{ccc}
\partial_{x_1}v_1 &   {1\over 2}(\partial_{x_1}v_2 + \partial_{x_2}v_1) &   {1\over 2}(\partial_{x_3}v_1 + \partial_{x_1}v_3)\\
\noame
 {1\over 2}(\partial_{x_1}v_2 + \partial_{x_2}v_1) & \partial_{x_2}v_2 &   {1\over 2}(\partial_{x_3}v_2 + \partial_{x_2}v_3)\\
 \noame
 {1\over 2}(\partial_{x_3}v_1 + \partial_{x_1}v_3)&   {1\over 2}(\partial_{x_3}v_2 + \partial_{x_2}v_3)& \partial_{x_3}v_3
\end{array}\right).$$

\noindent Moreover, for $\widetilde {  v} =(\widetilde {  v}', \widetilde  v_{3})$  a vector function and $\widetilde \varphi$ a scalar function, both defined in $\Omega$, obtained from ${ v}$ and $\varphi$ after a dilatation in the vertical variable, respectively, we will use the following operators
 $$\begin{array}{c}
 \displaystyle \Delta_{\varepsilon}  \widetilde  { v} =\Delta_{y'} \widetilde  { v} + \varepsilon^{-2}\partial_{y_3}^2 \widetilde  { v} ,\quad   \displaystyle\Delta_{ \varepsilon}\widetilde \varphi=\Delta_{y'}\widetilde \varphi+ \ep^{-2}\partial^2_{y_3}\widetilde \varphi,\\
  \noame
  (D_{\ep}\widetilde { v})_{ij}=\partial_{y_j}\widetilde {v}_i\ \hbox{ for }\ i=1,2,3,\ j=1,2,\quad   (D_{\ep}\widetilde { v})_{i3}=\ep^{-1}\partial_{y_3}\widetilde {v}_i \ \hbox{ for }\ i=1,2,3,\\
  \noame
  \nabla_{\ep}\widetilde\varphi=(\nabla_{y'}\widetilde \varphi, \ep^{-1}\partial_{y_3}\widetilde\varphi)^t,\quad {\rm div}_{\varepsilon}(\widetilde  {   v})={\rm div}_{y'}(\widetilde { v}')+\ep^{-1}{\partial_{y_3}}\widetilde  v_{3},
  \end{array}$$
 Moreover, we define $\mathbb{D}_{\ep}[\widetilde  { v}]$ as follows
 $$\mathbb{D}_{\ep}[\widetilde  { v}]=\mathbb{D}_{y'}[\widetilde  { v}]+\ep^{-1}\partial_{y_3}[\widetilde  { v}]=\left(\begin{array}{ccc}
\partial_{y_1}\widetilde  v_1 &   {1\over 2}(\partial_{y_1}\widetilde  v_2 + \partial_{y_2}\widetilde  v_1) &   {1\over 2} (\partial_{y_1}\widetilde  v_3+ \ep^{-1}\partial_{y_3}\widetilde  v_1)\\
\noame
 {1\over 2}(\partial_{y_1}\widetilde  v_2 + \partial_{y_2}\widetilde  v_1) & \partial_{x_2}\widetilde  v_2 &   {1\over 2} (\partial_{y_2}\widetilde  v_3+ \ep^{-1}\partial_{y_3}\widetilde  v_2)\\
 \noame
 {1\over 2} (\partial_{y_1}\widetilde  v_3+\ep^{-1}\partial_{y_3}\widetilde  v_1)&   {1\over 2} (\partial_{y_2}\widetilde  v_3+\ep^{-1}\partial_{y_3}\widetilde  v_2)& \ep^{-1}\partial_{y_3}\widetilde  v_3
\end{array}\right),$$
 where $\mathbb{D}_{y'}[\widetilde  { v}]$ and $\partial_{y_3}[\widetilde  { v}]$ are defined by
\begin{equation}\label{def_der_sym_1}
  \mathbb{D}_{y'}[{\widetilde v}]=\left(\begin{array}{ccc}
\partial_{y_1}\widetilde v_1 &   {1\over 2}(\partial_{y_1}\widetilde v_2 + \partial_{y_2}\widetilde v_1) &   {1\over 2} \partial_{y_1}\widetilde v_3\\
\noame
 {1\over 2}(\partial_{y_1}\widetilde v_2 + \partial_{y_2}\widetilde v_1) & \partial_{y_2}\widetilde v_2 &   {1\over 2} \partial_{y_2}\widetilde v_3\\
 \noame
 {1\over 2} \partial_{y_1}\widetilde v_3&   {1\over 2} \partial_{y_2}\widetilde v_3& 0
\end{array}\right),
 \  \partial_{y_3}[{ \widetilde v}]=\left(\begin{array}{ccc}
0 &   0&   {1\over 2} \partial_{y_3}\widetilde v_1\\
\noame
 0& 0 &   {1\over 2} \partial_{y_3}\widetilde v_2\\
 \noame
 {1\over 2} \partial_{y_3}\widetilde v_1&   {1\over 2} \partial_{y_3}\widetilde v_2& \partial_{y_3}\widetilde v_3
\end{array}\right).
\end{equation}

\noindent We also define the following operators applied to ${ v}'$:
\begin{equation}\label{def_der_sym_2}
  \mathbb{D}_{y'}[{\widetilde  v}']=\left(\begin{array}{ccc}
\partial_{y_1}\widetilde v_1 &   {1\over 2}(\partial_{y_1}\widetilde v_2 + \partial_{y_2}\widetilde v_1) &   0\\
\noame
 {1\over 2}(\partial_{y_1}\widetilde v_2 + \partial_{y_2}\widetilde v_1) & \partial_{y_2}\widetilde v_2 &  0\\
 \noame
0&   0& 0
\end{array}\right),
 \quad \partial_{y_3}[{\widetilde  v}']=\left(\begin{array}{ccc}
0 &   0&   {1\over 2} \partial_{y_3}\widetilde v_1\\
\noame
 0& 0 &   {1\over 2} \partial_{y_3}\widetilde v_2\\
 \noame
 {1\over 2} \partial_{y_3}\widetilde v_1&   {1\over 2} \partial_{y_3}\widetilde v_2& 0
\end{array}\right).
\end{equation}
We denote by $:$ the full contraction of two matrices, i.e. for $A=(a_{ij})_{1\leq i,j\leq 3}$ and $B=(b_{ij})_{1\leq i,j\leq 3}$, we have $A:B=\sum_{i,j=1}^3a_{ij}b_{ij}$.  The, we remark that 
$$ \mathbb{D}_{y'}[\widetilde w]: \mathbb{D}_{y'}[{\widetilde  v}']= \mathbb{D}_{y'}[{ w}']: \mathbb{D}_{y'}[{\widetilde  v}'],\quad \partial_{y_3}[{\widetilde  w}]:\partial_{y_3}[{\widetilde  v}']=\partial_{y_3}[{ \widetilde w}']:\partial_{y_3}[{\widetilde  v}'].$$

\noindent We denote by $C$ a generic constant which can change from line to line. Moreover, $O_\ep$ denotes a  generic quantity, which can change from line to line, devoted to tend to zero when $\ep\to 0$.

\section{Model problem}\label{sec:model} We assume that the the flow of velocity ${ u}_\ep=({ u}'_\ep(x), u_{3,\ep}(x))$ and pressure $p_\ep=p_\ep(x),$ at a point $x\in \Omega_\ep$, is suppose to be governed by the following nondimensional Navier-Stokes system (see e.g. \cite{Tapiero2, Tapiero})
\begin{equation}\label{system_1_NS}
\left\{\begin{array}{rl}
\displaystyle - {1\over Re} \,{\rm div}(\eta_r(\mathbb{D}[{u}_\ep])\mathbb{D}[{u}_\ep])+(u_\ep\cdot \nabla)u_\ep+\nabla p_\ep={f} & \hbox{in}\ \Omega_\varepsilon,\\
\noame
{\rm div}( {u}^\varepsilon)=0& \hbox{in}\ \Omega_\varepsilon\\
\noame
\displaystyle {u}_\varepsilon=0 & \hbox{on }\partial\Omega_\varepsilon,
\end{array}\right.
\end{equation}
where $\eta_r$ is the Carreau law without high rate viscosity defined  by (\ref{Carreaulaw}), and   in the momentum equation (\ref{system_1_NS})$_1$, the {\it Reynolds number} $Re$ is assumed to be proportional to $\varepsilon^{-1}$, i.e. $Re=\varepsilon^{-1}$.  Moreover, following \cite{Tapiero2}, the source term ${f}$ is of the form
\begin{equation}\label{fassump}
{f} (x)=({f}'(x'),0)\quad \hbox{with}\quad {f}'\in L^\infty(\omega)^2.
\end{equation}
Before studying the limit behavior of the above system when $\ep\to 0$, we discuss solvability of system (\ref{system_1_NS}) with $\eta_\infty=0$, which can be found in \cite{Mikelicmu0}.  For every $\ep>0$, it is well known that this problem has at least one variational weak solution $(u_\ep, p_\ep)$ in $W^{1,r}_0(\Omega_\ep)\times L^{r'}_0(\Omega_\ep)$ for $ 9/5<r<+\infty$ (see \cite{Lions2}), which was latter extended to $6/5<r<+\infty$ (see \cite{FMS}),  where  $L^{r'}_0$  is the space of functions of $L^{r'}$   with zero mean value. 
\begin{remark}We recall the difference with the the case $\eta_\infty>0$ and $1<r<+\infty$, $r\neq 2$,  in which  there exists of at least one weak solution of system (\ref{system_1_NS}) given by $(u_\varepsilon, p_\epsilon) \in H^1_0(\Omega_\varepsilon)^3\times L^2_0(\Omega_\epsilon)$, for $1<r<2$,  and $(u_\varepsilon, p_\epsilon) \in W^{1,r}_0(\Omega_\epsilon)^3\times L^{r'}_0(\Omega_\varepsilon)$ wtih $1/r+1/r'=1$, for $r>2$.
\end{remark}

Then, to study the limit behavior as $\ep\to 0$,  the usual procedure is to compute the a priori estimates after testing equation (\ref{system_1_NS}) with the solution, which gives (see  Proposition \ref{estim_sol_dil}  below) 
$$\|u_\ep\|_{L^r(\Omega_\ep)^3}+\ep\|D u_\ep\|_{L^r(\Omega_\ep)^{3\times 3}}\leq C\ep^{1+{1\over r}},$$
i.e. $u_\ep$ is of order $\ep^{1+{1\over r}}$. As a consequence, it can be shown that, when we study the limiting behavior as $\ep\to 0$ after a dilatation in the vertical variable, the inertial term vanishes in the limit. So it is justified to consider the quasi-Newtonian Stokes system (for similar  reductions  from Navier-Stokes to Stokes,  see e.g. \cite{Gotz, MikelicIntro2}). Consequently, from now on, we consider the dimensionless system
\begin{equation}\label{system_1}
\left\{\begin{array}{rl}
\displaystyle - \varepsilon \,{\rm div}(\eta_r(\mathbb{D}[{u}_\ep])\mathbb{D}[{u}_\ep])+\nabla p_\ep={f} & \hbox{in}\ \Omega_\varepsilon,\\
\noame
{\rm div}( {u}^\varepsilon)=0& \hbox{in}\ \Omega_\varepsilon,\\
\noame
\displaystyle {u}_\varepsilon=0 & \hbox{on }\partial\Omega_\varepsilon.
\end{array}\right.
\end{equation}
We remark, that for the Stokes system, for every $\ep>0$, the classical method of monotone operators, see for instance \cite{Baranger,  BourgeatPaloka0, Sandri} in the case $\eta_\infty=0$,  gives the existence and uniqueness of the velocity $u_\varepsilon\in W^{1,r}_0(\Omega_\varepsilon)^3$ and of the pressure $p_\varepsilon \in L^{r'}_0(\Omega_\varepsilon)$. 
\begin{remark} For Stokes system (\ref{system_1}),  we recall the difference with the the case $\eta_\infty>0$ and $1<r<+\infty$, $r\neq 2$,  in which  there exists a unique weak  solution $u_\varepsilon\in H^1_0(\Omega_\varepsilon)^3$ and  $p_\varepsilon\in L^{2}_0(\Omega_\varepsilon)$, see for instance \cite{Baranger,  Tapiero2, BourgeatPaloka0, Sandri}.
\end{remark}
To study the asymptotic behavior of the solutions ${ u}_\ep$ and $p_\ep$ when $\ep$  tends to zero,  we use the rescaling  
\begin{equation}\label{dilatacion}
y'=x',\quad y_3={x_3\over \ep}\,,
\end{equation}
 to have the functions defined in  $\Omega$, which is defined in (\ref{domains_tilde}). The rescaled  Stokes system is given as follows
\begin{equation}\label{system_1_dil}
\left\{\begin{array}{rl}
\displaystyle -   \varepsilon \,{\rm div}_\varepsilon(\eta_r(\mathbb{D}_\varepsilon[ \widetilde { u}_\ep])\mathbb{D}_\varepsilon[\widetilde {u}_\ep])+\nabla_\varepsilon \widetilde p_\ep={f} & \hbox{in}\   \Omega,\\
\noame
{\rm div}_\varepsilon( \widetilde {u}_\varepsilon)=0& \hbox{in}\  \Omega,\\
\noame
\displaystyle \widetilde {u}_\varepsilon=0 & \hbox{on }\partial \Omega,
\end{array}\right.
\end{equation}

\noindent where the unknown functions in the above system are given by 
$${\widetilde u}_\varepsilon(y)={ u}_\varepsilon(y',\varepsilon y_3),\quad \widetilde p_\varepsilon(y)=p_\varepsilon(y',\varepsilon y_3),\quad  \hbox{for a.e. }y\in  \Omega,$$ and the operators ${\rm div}_\ep$, $\mathbb{D}_\ep$ and $\nabla_\ep$ are defined in Section \ref{sec:definition}.

\section{The formal asymptotic expansion}\label{sec:formal}

In this section, we apply the asymptotic expansion method (see \cite{Bayada_Chambat_Transition, MikelicIntro2} for instance) to the problem (\ref{system_1_dil}) to derive a reduced Stokes model. The idea is to assume an expansion in $\ep$ of the solution $(\widetilde u_\ep, \widetilde p_\ep)$ given by
\begin{equation}\label{expansion}
\widetilde u_\ep(y)=\ep^\beta\left(v^0(y)+\ep v^1(y)+\ep^2v^2(y)+\cdots\right),\quad \widetilde p_\ep(y)= p^0(y)+\ep p^1(y)+\ep^2p^2(y)+\cdots
\end{equation}
To determine the limit problems given by functions $(v^0,p^0)$, the expansion (\ref{expansion}) is plugged into the PDE, we identify the various powers of $\ep$ and we obtain a cascade of equations from which we retain only the leading ones that constitute the limit problem.

We remark that the value of $\beta$ in the expansion (\ref{expansion}) will be determined in the next section, by the {\it a priori} estimates for $(\widetilde u_\ep, \widetilde p_\ep)$ (see Lemmas \ref{lemma_estimates} and \ref{lemma_est_P}). Moreover, the limit problems will be justified by corresponding compactness results (see Lemma \ref{lem_conv_vel} and Theorem \ref{them_limit}).

Thus, consider $(\widetilde u_\ep, \widetilde p_\ep)$ the solution of  problem (\ref{system_1_dil}). Assuming the asymptotic expansion of the unknown $(\widetilde u_\ep, \widetilde p_\ep)$ in the following form 
\begin{equation}\label{expansion_thm}
\widetilde u_\ep(y)=\ep v^0(y)+\ep^2v^1(y)+\cdots,\quad p_\ep(y)= p^0(y)+\ep p^1(y)++\cdots \quad\hbox{for a.e. }y\in \Omega,
\end{equation}
where  $v^i(y)=(\overline v^i(y), v_3^i(y))$, with $\overline v^i(y)=(v^i_1(y),v^i_2(y))$, $i=0,1,\ldots$, we obtain:\\

\begin{itemize}
\item[--]  {\it Boundary conditions}. By using expansion (\ref{expansion_thm}), we first derive boundary conditions on $\Gamma_0\cup \Gamma_1$. The main order and next order terms in the boundary condition $\widetilde u_\ep=0$ on $\Gamma_0\cup\Gamma_1$ are
\begin{equation}\label{bc_expansion}
\begin{array}{lll}
\ep^{1}:& v^0=0&\hbox{on }\Gamma_0\cup\Gamma_1,\\
\noame
\ep^{2}:& v^1=0&\hbox{on }\Gamma_0\cup\Gamma_1,\\
\end{array}
\end{equation}
From the previous equalities, we deduce the boundary conditions 
 $$\overline v^0=0\quad\hbox{on }\Gamma_0\cup \Gamma_1.$$
 
\item[--]  {\it Divergence condition.} Next,  the expansion in the divergence condition ${\rm div}_\ep(\widetilde u_\ep)=0$ in $\Omega$ gives
$$\ep {\rm div}_{y'}(\overline v^0+\ep \overline v^1+\cdots )+ \partial_{y_3}( v^0_3+\ep  v^1_3+\cdots )=0,\quad \hbox{in }\Omega,$$
and so,  the main order terms are the following:
\begin{equation}\label{div_expansion}
\begin{array}{lll}
1:& \partial_{y_3}  v^0_3=0&\hbox{in }\Omega,\\
\noame
\ep :& {\rm div}_{y'}(\overline v^0)+\partial_{y_3}v_3^1=0&\hbox{in }\Omega.\\
\end{array}
\end{equation}
According to  (\ref{div_expansion})$_1$, $v_3^0$ does not depend on $y_3$, and the boundary conditions for $v_3^0$ given in (\ref{bc_expansion})$_1$, we deduce that $v_3^0\equiv 0$. Also, integrating (\ref{div_expansion})$_2$ with respect to $y_3$ between $0$ and $h(y')$ and using the boundary conditions of $v_3^1$ given in (\ref{bc_expansion})$_2$, we deduce the incompressibility condition 
$${\rm div}_{y'}\left(\int_0^{h(y')}\overline v^0\,dy_3\right)=0\quad\hbox{in }\omega.$$

\item[--] {\it Equation for $\overline v^0$.} Now, substituting the expansion into the equation (\ref{system_1_dil})$_1$ and taking into account $v^0_3\equiv 0$, we get for the horizontal variables
\begin{equation}\label{mom_expansion2}
\begin{array}{l}\displaystyle
-\eta_0\ep  {\rm div}_{y'}\left(\left(1+\lambda |\ep\mathbb{D}_\ep[ \overline v^0+O(\ep)]|^2\right)^{r-2\over 2}(\ep\mathbb{D}_\ep[\overline v^0+O(\ep))]\right)\\
\noame
\displaystyle
- \eta_0  \partial_{y_3}\left(\left(1+\lambda |\ep \mathbb{D}_\ep[ \overline v^0+O(\ep)]|^2\right)^{r-2\over 2}(\ep\mathbb{D}_\ep[\overline   v^0+O(\ep)])
\right)\\
\noame
\displaystyle
+\nabla_{y'}(p^0+O(\ep)) =f',
\end{array}
\end{equation}
and the equation for the vertical variable is the following
\begin{equation}\label{mom_expansion_3}
\begin{array}{l}\displaystyle
-\eta_0\ep^2 {\rm div}_{y'}\left(\left(1+\lambda|\ep \mathbb{D}_\ep[ \overline v^0+O(\ep)]|^2\right)^{r-2\over 2}(\ep \mathbb{D}_\ep[ v_3^1+ O(\ep)])\right)\\
\noame
\displaystyle
- \eta_0\ep  \partial_{y_3}\left(\left(1+\lambda|\ep \mathbb{D}_\ep[ \overline v^0+O(\ep)]|^2\right)^{r-2\over 2}(\ep \mathbb{D}_\ep[  v^1_3+O(\ep)])
\right)\\
\noame
\displaystyle
 +\ep^{-1}\partial_{y_3}(p^0+\ep p^1 + O(\ep^2))=0.
\end{array}
\end{equation}
From (\ref{mom_expansion2}) and (\ref{mom_expansion_3}), since $p^0$ does not depend on $y_3$, we get
$$
\begin{array}{lll}
1:& \partial_{y_3}  p^1=0&\hbox{in }\Omega,\\
\noame
1:& -\eta_0  \partial_{y_3}\left(\left(1+\lambda |\partial_{y_3}[\overline v^0]|^2\right)^{r-2\over 2}\partial_{y_3}[\overline v^0]
\right) 
+\nabla_{y'}p^0 =f'&\hbox{in }\Omega,
\end{array}
$$
which gives 
$$-\partial_{y_3}\left(\left(1+{\lambda\over 2}|\partial_{y_3} \overline v^0|^2\right)^{{r\over 2}-1} \partial_{y_3}\overline v^0\right)={2\over  {\eta_0} }({f}'(y')-\nabla_{y'}  p^0(y')) \quad\hbox{in }\Omega.$$
\end{itemize}
In summary, we deduce that the main order pair of functions $(v^0,p^0)$ with $v_3^0\equiv 0$ and $p^0=p^0(y')$, satisfies the following reduced Carreau Stokes model
\begin{equation}\label{limit_model_formal}
\left\{\begin{array}{rl}
\displaystyle
-\partial_{y_3}\left(\left(1+{\lambda\over 2}|\partial_{y_3} \overline v^0|^2\right)^{{r\over 2}-1} \partial_{y_3}\overline v^0\right)={2\over  {\eta_0} }({f}'(y')-\nabla_{y'}  p^0(y')) &\hbox{in }\Omega,\\
\noame
\displaystyle
{\rm div}_{y'}\left(\int_0^{h(y')}\overline v^0\,dy_3\right)=0&\hbox{in }\omega,\\
\noame
\displaystyle  \overline v^0=0&\hbox{on }\Gamma_0\cup \Gamma_1.
\end{array}\right.
\end{equation}

In the next section, we derive the sharp {\it a priori} estimates for velocity and pressure, which will allow us to later prove  compactness results of the rescaled functions.

\section{{\it A priori} estimates} \label{sec:estimates}

 Let us recall and proof some well-known technical results (see for instance \cite[Lemmas 0.2 and 0.3]{Tapiero2}).
\begin{lemma}[Poincar\'e's  and Korn's inequalities]\label{Poincare_lemma} For every $\varphi\in W^{1,q}_0(\Omega_\varepsilon)^3$, $1< q<+\infty$, it holds
\begin{equation}\label{Poincare}
\|\varphi\|_{L^q(\Omega_\varepsilon)^3}\leq C\varepsilon\|D \varphi\|_{L^q(\Omega_\varepsilon)^{3\times 3}},
\quad
 \|D\varphi\|_{L^q(\Omega_\varepsilon)^{3\times 3}}\leq C\|\mathbb{D}[\varphi]\|_{L^q(\Omega_\varepsilon)^{3\times 3}}.
\end{equation}
From the change of variables (\ref{dilatacion}),    for every $\widetilde\varphi\in W^{1,q}_0(\Omega)^3$, it holds  
\begin{equation}\label{Poincare2}
\|\widetilde \varphi\|_{L^q(\Omega)^3}\leq C\varepsilon\|D_{\ep} \widetilde \varphi\|_{L^q(\Omega)^{3\times 3}},\quad
 \|D_\varepsilon\widetilde \varphi\|_{L^q(\Omega)^{3\times 3}}\leq C\|\mathbb{D}_\varepsilon[\widetilde \varphi]\|_{L^q(\Omega)^{3\times 3}}.
\end{equation}
\end{lemma}

\begin{lemma}[Velocity estimates]\label{lemma_estimates} We have the following estimates for the solution ${\widetilde u}_\ep$ of  (\ref{system_1_dil}):
\begin{equation}\label{estim_sol_dil}
\displaystyle
\|{\widetilde u}_\varepsilon\|_{L^r(\Omega)^3}\leq C\ep, \quad\displaystyle
\|D_{\varepsilon} { \widetilde u}_\varepsilon\|_{L^r(\Omega)^{3\times 3}}\leq C,\quad\displaystyle
\|\mathbb{D}_{\varepsilon} [{ \widetilde u}_\varepsilon]\|_{L^r(\Omega)^{3\times 3}}\leq C.
\end{equation}
\end{lemma}
\begin{proof}
Multiplying (\ref{system_1_dil})$_1$ by $\ep^{-1}\widetilde u_\varepsilon$, integrating over $\Omega_\ep$ and taking into account  that ${\rm div}(u_\varepsilon)=0$ in $\Omega_\ep$, the Poincar\'e and Korn inequalities (\ref{Poincare}) and assumptions on $f$,  we get 
\begin{equation}\label{estim_1_dev} \eta_0 \int_{\Omega}(1+\lambda|\mathbb{D}_\ep[{\widetilde u}_\ep]|^2)^{{r\over 2}-1}|\mathbb{D}_\ep [{\widetilde u}_\ep]|^2\,dx\leq C\ep^{-1}\|\widetilde  u_\ep\|_{L^r(\Omega)^{3}}\leq C \|\mathbb{D}_\ep[\widetilde u_\ep]\|_{L^r(\Omega)^{3\times 3}}.
\end{equation}
Consider the case $1<r<2$. After the application of Jensen's inequality to the convex function $\Phi(z)={z^{2\over r}\over 1+(\lambda z)^{2-r\over r}}$ (see \cite[Proposition 3.1]{Baranger} or \cite[Proposition 2.4]{Mikelicmu0} for instance) we obtain
$$  \begin{array}{rl}\displaystyle \int_{\Omega}(1+\lambda|\mathbb{D}_\ep[{\widetilde u}_\ep]|^2)^{{r\over 2}-1}|\mathbb{D}_\ep [{\widetilde u}_\ep]|^2\,dx\geq &\displaystyle \int_{\Omega}{|\mathbb{D}_\ep[{\widetilde u}_\ep]|^2\over 1+\lambda^{2-r\over 2}|\mathbb{D}_\ep[{\widetilde u}_\ep]|^{2-r}}\geq\displaystyle\alpha_1{\|\mathbb{D}_\ep[\widetilde u_\ep]\|_{L^r(\Omega)^{3\times 3}}^2\over 1+\alpha_2\|\mathbb{D}_\ep[\widetilde u_\ep]\|_{L^r(\Omega)^{3\times 3}}^{2-r}},
\end{array}$$
with $\alpha_1, \alpha_2>0$. This with (\ref{estim_1_dev}) implies
$$\|\mathbb{D}_\ep[\widetilde u_\ep]\|_{L^r(\Omega)^{3\times 3}}^2\leq C\|\mathbb{D}_\ep[\widetilde u_\ep]\|_{L^r(\Omega)^{3\times 3}}(1+\|\mathbb{D}_\ep[\widetilde u_\ep]\|_{L^r(\Omega)^{3\times 3}}^{2-r}),$$
and so, it holds
$$\|\mathbb{D}_\ep[\widetilde u_\ep]\|_{L^r(\Omega)^{3\times 3}}\leq C(1+\|\mathbb{D}_\ep[\widetilde u_\ep]\|_{L^r(\Omega)^{3\times 3}}^{2-r}).$$
By using Young's inequality with $q=1/(2-r)$ and $q'=1/(r-1)$ to $C\|\mathbb{D}_\ep[\widetilde u_\ep]\|_{L^r(\Omega)^{3\times 3}}^{2-r}$, we deduce 
$$\|\mathbb{D}_\ep[\widetilde u_\ep]\|_{L^r(\Omega)^{3\times 3}}\leq C+{C^{q'}\over q'}\left({2\over q}\right)^{q'\over q}  +{1\over 2}\|\mathbb{D}_\ep[\widetilde u_\ep]\|_{L^r(\Omega)^{3\times 3}},$$
which gives (\ref{estim_sol_dil})$_{3}$.

Now, we consider the case $r>2$. Due to the inequality
$$  \begin{array}{rl}\displaystyle \int_{\Omega}(1+\lambda|\mathbb{D}_\ep[{\widetilde u}_\ep]|^2)^{{r\over 2}-1}|\mathbb{D}_\ep [{\widetilde u}_\ep]|^2\,dx\geq &\displaystyle  \lambda^{r-2\over 2}\int_{\Omega} |\mathbb{D}_\ep[{\widetilde u}_\ep]|^{r}\,dx,
\end{array}$$
and (\ref{estim_1_dev}), we deduce 
$$\|\mathbb{D}_\ep[\widetilde u_\ep]\|_{L^r(\Omega)^{3\times 3}}^r\leq C\|\mathbb{D}_\ep[\widetilde u_\ep]\|_{L^r(\Omega)^{3\times 3}},$$
which implies (\ref{estim_sol_dil})$_{3}$.
  
  Finally, from the Poincar\'e and Korn inequelities (\ref{Poincare}) applied to (\ref{estim_sol_dil})$_3$, we deduce  (\ref{estim_sol_dil})$_{1,2}$.

\end{proof}

Next, we decompose the pressure $p_\varepsilon$, by using a result given in \cite{CLS}, in two pressures $p_\ep^0$ and $p_\ep^1$ and we derive  the corresponding estimates. As a result, we have that the main pressure $p_\ep^0$ does not depend on $y_3$ and have regularity in $W^{1,r'}$ instead of $L^{r'}$, which will allow to have strong convergence of the pressure.
\begin{lemma}\label{lemma_est_P} The pressure solution $p_\ep\in L^r_0(\Omega_\ep)$ of (\ref{system_1})  can be decomposed as follows
\begin{equation}\label{decomposition_original}p_\ep=  p_\ep^0+  p_\ep^1,
\end{equation}
with $p_\ep^0\in W^{1,r'}(\omega)$ and $p_\ep^1\in L^{r'}(\Omega_\ep)$ with
\begin{equation}\label{estim_P_original2}
\|\nabla_{x'}p_\ep^0\|_{L^{r'}(\omega)}\leq C,\quad \|p_\ep^1\|_{L^{r'}(\Omega_\ep)}\leq C\ep^{{1\over r'}+1}.
\end{equation}
Moreover, by applying the change of variables (\ref{dilatacion}), we have  
\begin{equation}\label{esti_P}
 \|\widetilde p_\ep^1\|_{L^{r'}(\Omega)^3}\leq C\ep.
\end{equation}
\end{lemma}

\begin{proof} 

We divide the proof in two steps. In the first step, we decompose the pressure $p_\ep$ of the pressure in the sum of two different pressures $p_\ep^0$ and $p_\ep^1$ and give estimates of both pressures with respect to an estimate of $\nabla p_\ep$ in $W^{-1,r'}(\Omega_\ep)^3$. In the second step, we derive estimates for $\nabla p_\ep$, and as consequence, for $p_\ep^0$ and $p_\ep^1$.

{\it Step 1. Decomposition of the pressure.} According to   \cite{CLS}, the pressure $p_\ep\in L^{r'}_0(\Omega_\ep)$ can be decomposed as in (\ref{decomposition_original}) with $p_\ep^0\in W^{1,r'}(\omega)$ and $p_\ep^1\in L^{r'}(\Omega_\ep)$, satisfying the following estimates 
\begin{equation}\label{estim_decomposition}\ep^{{1\over r'}+1}\|\nabla_{x'}p_\ep^0\|_{L^{r'}(\omega_\ep)}+\|p_\ep^1\|_{L^{r'}(\Omega_\ep)}\leq \|\nabla p_\ep\|_{W^{-1,r'}(\Omega_\ep)}.
\end{equation}
 We point out that in \cite[Corollary 3.4]{CLS}, the scaling of the pressure is $p_\ep={1\over \ep}\pi_\ep^0+\pi_\ep^1$ satisfying
$$\ep^{{1\over r'}}\|\nabla_{x'}\pi_\ep^0\|_{L^{r'}(\omega)}+  \|\pi_\ep^1\|_{L^{r'}(\Omega_\ep)}\leq C\|\nabla p_\ep\|_{W^{-1,r'}(\Omega_\ep)^3}.$$
Here, we rescale  $p_\ep^0=\ep^{-1} \pi_\ep^0$ and $p_\ep^1= \pi_\ep^1$ and so, we get (\ref{estim_decomposition}).
 \\

{\it Step 2. Estimates for $p_\ep^0$ and $p_\ep^1$.} We derive estimates for pressure by using the variational formulation of problem (\ref{system_1}), i.e. equality
  \begin{equation}\label{final_step1}\langle\nabla  p_\varepsilon, \varphi\rangle_{W^{-1,r'}(\Omega_\ep)^3,W^{1,r}_0(\Omega_\ep)^3}=- \ep \eta_0\int_{\Omega_\ep}(1+\lambda|\mathbb{D}[{ u}_\ep]|^2)^{{r\over 2}-1}\mathbb{D}[{ u}_\ep]:\mathbb{D}[\varphi]\,dx+\int_{\Omega_\ep} {f}'\cdot   \varphi'\,dx,
  \end{equation}
  for every $\varphi\in W^{1,r}_0(\Omega_\ep)^3$.  Using estimates for velocity $u_\ep$, we deduce
\begin{equation}\label{Step2_2} 
\displaystyle \left|\ep  \eta_0 \int_{\Omega_\ep}(1+\lambda|\mathbb{D}[{ u}_\ep]|^2)^{{r\over 2}-1}\mathbb{D} [{ u}_\ep]:\mathbb{D}[\varphi]\,dx\right| \leq C\ep^{{1\over r'}+1}\|D\varphi\|_{L^r(\Omega_\ep)^{3\times 3}}\leq C\ep^{{1\over r'}+1}\|\varphi\|_{W^{1,r}_0(\Omega_\ep)^3}.
\end{equation}

Since ${ f}'={ f}'(x')$ is in $L^\infty(\omega)^2$ and using Poincare's inequality (\ref{Poincare})$_1$, we get
\begin{equation}\label{Step2_5} 
\displaystyle \left| \int_{\Omega_\ep} {f}'\cdot   \varphi'\,dy\right| \leq   C\ep^{1\over r'}\|\varphi\|_{L^{r}(\Omega_\ep)^3}  \leq   C\ep^{{1\over r'}+1}\|D\varphi\|_{L^{r}(\Omega_\ep)^3} \leq C\ep^{{1\over r'}+1}\|\varphi\|_{W^{1,r}_0(\Omega)^3}.
\end{equation}
Coming back to the expression (\ref{final_step1}), we deduce from (\ref{Step2_2})-(\ref{Step2_5}) the estimate
$$\|\nabla p_\ep\|_{W^{-1,r'}(\Omega_\ep)^3}\leq C\ep^{{1\over r'}+1}.$$
This estimate and (\ref{estim_decomposition}) give estimates (\ref{estim_P_original2}).

Finally, by using the change of variables (\ref{dilatacion}), we deduce estimate (\ref{esti_P}) from (\ref{estim_P_original2})$_2$.

\end{proof}

\section{Convergences and limit models}\label{sec:convergences}
In this section, we analyze the asymptotic behavior of functions $(\widetilde { u}_\ep, \widetilde p_\ep)$ when $\ep$ tends to zero. We introduce the space $V_{y_3}^r(\Omega)$   defined by
$$V_{y_3}^r(\Omega)=\{\varphi\in L^r(\Omega)\ :\ \partial_{y_3}\varphi\in L^r(\Omega)\},$$
which is a Banach space with the norm (see for instance \cite{Boukrouche_ElMir})
$$\varphi\to \|\varphi\|_{V_{y_3}^r(\Omega)}=\left(\|\varphi\|_{L^r(\Omega)^3}+\|\partial_{y_3}\varphi\|_{L^r(\Omega)^3}
\right)^{1\over r}.$$

\begin{lemma} \label{lem_conv_vel}There exist $\widetilde { u}\in V_{y_3}^r(\Omega)^3$ where $\widetilde { u}=0$ on $\Gamma_0\cup \Gamma_1$ and $\widetilde u_3\equiv 0$, and $\widetilde p\in L^{r'}_0(\omega)\cap W^{1,r'}(\omega)$, independent of $y_3$ such that, up to a subsequence, 
\begin{equation}\label{conv_vel_tilde}
\ep^{-1}\widetilde { u}_\ep\rightharpoonup (\widetilde { u}',0)\quad\hbox{in }V_{y_3}^r(\Omega)^3,
\end{equation}
\begin{equation}\label{pressuresr01}
 p_\ep^0\rightharpoonup \widetilde p\quad\hbox{in }W^{1,r'}(\omega),\quad \ep^{-1}\widetilde p_\ep^1\rightharpoonup \widetilde p_1\quad\hbox{in }L^{r'}(\Omega),
\end{equation}
\begin{equation}\label{pressuresr}
\widetilde p_\ep\to \widetilde p\quad\hbox{in }L^{r'}(\Omega), 
\quad \ep^{-1}\partial_{y_3}\widetilde p_\ep \rightharpoonup \widetilde p_1\quad\hbox{in }L^{r'}(\omega;W^{-1,r'}(0,1)).
\end{equation}
\begin{equation}\label{div_conv}
{\rm div}_{y'}\left(\int_0^{h(y')}\widetilde { u}'(y)\,dy_3\right)=0\quad \hbox{in }\omega,\quad\left(\int_0^{h(y')}\widetilde { u}'(y)\,dy_3\right)\cdot n=0\quad \hbox{on }\partial\omega.
\end{equation}

\end{lemma}
 \begin{proof}
The proof (\ref{conv_vel_tilde}) and (\ref{div_conv}) is a direct consequence of estimates (\ref{estim_sol_dil}), see for instance \cite{Tapiero2, Tapiero}, so we omit it.

According to estimates (\ref{estim_P_original2})$_1$ and (\ref{esti_P}), we deduce that there exists $\widetilde p\in W^{1,r'}(\omega)$ and $\widetilde p_1\in L^{r'}(\Omega)$ such that convergences given in (\ref{pressuresr01})  hold. As consequence, it holds  convergences  (\ref{pressuresr}).

Finally, since $\widetilde p_\ep\in L^{r'}_0(\Omega)$, then
$$0=\int_{\Omega}\widetilde p_\ep\,dy=\int_{\Omega}  p_\ep^0\,dy+\int_{\Omega}\widetilde p_\ep^1\,dy,$$
and by using convergences (\ref{pressuresr01}), we deduce
$$0=\int_{\Omega}\widetilde p\,dy=h(y')\int_{\omega}\widetilde p\,dy',$$
so $\widetilde p\in L^{r'}_0(\omega)$.

 \end{proof}

Next, using monotonicity arguments together with Minty's lemma (see for instance \cite{Tapiero2, EkelandTemam}), we  derive a variational inequality that will be useful in the proof of the main theorem.  
 
According to Lemma \ref{lem_conv_vel}, we choose a test function ${ v}(y)\in \mathcal{D}(\Omega)^3$ with $v_3\equiv 0$. Multiplying (\ref{system_1_dil}) by ${ v}(y)$, using the decomposition of the pressure and integrating by parts, we have
\begin{equation}\label{form_var_hat}\begin{array}{l}
\displaystyle
\medskip
  \eta_0\int_{ \Omega }(1+\lambda|\mathbb{D}_\varepsilon[\widetilde { u}_\varepsilon]|^2)^{{r\over 2}-1}\mathbb{D}_\varepsilon[\widetilde {u}_\varepsilon] :(\ep \mathbb{D}_{y'}[ { v}]+\partial_{y_3}[v])\,dy
\\
\noame
\displaystyle -\int_{\Omega}p_\ep^0 {\rm div}_{y'}(v')\,dy-\int_{\Omega}\widetilde p_\varepsilon^1  {\rm div}_{\ep}({ v})\,dy=\int_{\Omega} { f}'\cdot { v}'\,dy.
\end{array}
\end{equation}

Now, let us define the functional $J_r$ by
$$
J_r({ v})=2{\eta_0\over r\lambda}\int_{\Omega}(1+\lambda| \varepsilon\mathbb{D}_{y'}[  { v}]+\partial_{y_3}[v] |^2)^{r\over 2}dy.
$$
Observe that, for every $\varepsilon>0$,  $J_r$ is convex and Gateaux differentiable on $W^{1,r}_0(\Omega)^3$, (see  \cite[Proposition 2.1 and Section 3]{Baranger}) and $A_r=J'_r$ is given by
$$\begin{array}{l}\displaystyle
(A_r({ w}),{ v})
=\displaystyle \eta_0\int_{\Omega}(1+\lambda| \ep \mathbb{D}_{y'}[  { w}]+ \partial_{y_3}[v]|^2)^{{r\over 2}-1}(\ep\mathbb{D}_{y'}[  { w}]+\partial_{y_3}[v]): (\ep \mathbb{D}_{y'}[  { v}] +\partial_{y_3}[v])\,dy.
\end{array}
$$
Applying \cite[Proposition 1.1., p.158]{Lions2}, $A_r$ is monotone, \emph{i.e.}
\begin{equation}\label{monotonicity}
(A_r({w})-A_r({ v}),{ w}-{ v})\ge0,\quad \forall { w},{ v}\in W^{1,r}_0(\Omega)^3.
\end{equation}
On the other hand, for all $\varphi\in \mathcal{D}(\Omega)^3$ with $\varphi_3\equiv 0$,  we choose ${ v}_\varepsilon$ defined by
\begin{equation}\label{testv}
{ v}_\varepsilon=\varphi-\varepsilon^{-1}\widetilde { u}_\varepsilon,
\end{equation}
as a test function in (\ref{form_var_hat}), and so we  have
\begin{equation*}\begin{array}{l}
\displaystyle
\medskip
 (A_r(\varepsilon^{-1}\widetilde { u}_\varepsilon),{ v}_\varepsilon)-\int_{\Omega}p_\ep^0 {\rm div}_{y'}({ v}'_\varepsilon)\,dy-\int_{\Omega}\widetilde p_\varepsilon^1  {\rm div}_\ep({ v}_\varepsilon)\,dy=\int_{\Omega} { f}'\cdot { v}'_\varepsilon\,dy,
\end{array}
\end{equation*}
which is equivalent to
\begin{equation*}\begin{array}{l}
\displaystyle
\medskip
 (A_r(\varphi)-A_r(\varepsilon^{-1}\widetilde {u}_\varepsilon),{ v}_\varepsilon)-   (A_r(\varphi),{ v}_\varepsilon)\\
 \noame
 \displaystyle
 +\int_{\Omega}p_\ep^0 {\rm div}_{y'}({ v}'_\varepsilon)\,dy+\int_{\Omega}\widetilde p_\varepsilon^1  {\rm div}_\ep({ v}_\varepsilon)\,dy=-\int_{\Omega} { f}'\cdot { v}'_\varepsilon\,dy.
\end{array}
\end{equation*}
Due to (\ref{monotonicity}), we can deduce
\begin{equation*}\begin{array}{l}
\displaystyle
\medskip
   (A_r(\varphi),{ v}_\varepsilon)-\int_{\Omega}p_\ep^0 {\rm div}_{y'}({ v}'_\varepsilon)\,dy-\int_{\Omega}\widetilde p_\varepsilon^1  {\rm div}_\ep({ v}_\varepsilon)\,dy\ge\int_{\Omega} {f}'\cdot {v}'_\varepsilon\,dy,
\end{array}
\end{equation*}
\emph{i.e.} we have  (recall $\varphi_3\equiv 0$)  that
\begin{equation}\label{v_ineq_carreau}\begin{array}{l}
\displaystyle
\medskip
   \eta_0\int_{\Omega}(1+\lambda|\ep\mathbb{D}_{y'}[  \varphi'] + \partial_{y_3}[ \varphi']|^2)^{{r\over 2}-1}(\varepsilon  \mathbb{D}_{y'}[  \varphi'] + \partial_{y_3}[ \varphi']):(\varepsilon \mathbb{D}_{y'}[  { v}_\ep'] + \partial_{y_3}[ { v}_\ep'])dy\\
\noame
\displaystyle-\int_{\Omega}p_\ep^0 {\rm div}_{y'}({ v}'_\varepsilon)\,dy-\int_{\Omega}\widetilde p_\varepsilon^1  {\rm div}_\ep({ v}_\varepsilon)\,dy\ge\int_{\Omega} { f}'\cdot {v}'_\varepsilon\,dy.
\end{array}
\end{equation}
Now, according to the compactness results given in Lemma \ref{lem_conv_vel}, we give the main result of the paper which gives the limit model.

\begin{theorem}[Limit model]\label{them_limit}  The sequence  $(\widetilde u_\ep, \widetilde p_\ep)$ of solution of (\ref{system_1_dil}) satisfies 
\begin{equation}\label{convergences_thm}\begin{array}{l}\displaystyle 
\ep \widetilde u_\ep\rightharpoonup 0 \quad\hbox{in } W^{1,r}(\Omega)^3,\quad \ep^{-1}\widetilde u_\ep\rightharpoonup \widetilde u \quad\hbox{in } V_{y_3}^r(\Omega)^3,\\
\noame
\displaystyle \widetilde p_\ep\to \widetilde p\quad\hbox{in }L^{r'}(\Omega),\quad \ep^{-1}\partial_{y_3}\widetilde p_\ep \rightharpoonup 0 \quad \hbox{in }L^{r'}(\omega;W^{-1,r'}(0,1)),
\end{array}
\end{equation}
where the pair $(\widetilde u, \widetilde p)\in V_{y_3}^r(\Omega)^3\times (L^{r'}_0(\omega)\cap W^{1,r'}(\omega))$, with $\widetilde u_3\equiv0$ satisfies  the following reduced Carreau Stokes system without high rate viscosity  
\begin{equation}\label{limit_model}
\left\{\begin{array}{rl}
\displaystyle
-\partial_{y_3}\left(\left(1+{\lambda\over 2}|\partial_{y_3} \widetilde { u}'|^2\right)^{{r\over 2}-1} \partial_{y_3}\widetilde {u}'\right)={2\over  {\eta_0} }({f}'(y')-\nabla_{y'}\widetilde p(y')) &\hbox{in }\Omega,\\
\noame
\displaystyle
{\rm div}_{y'}\left(\int_0^{h(y')}\widetilde u'\,dy_3\right)=0&\hbox{in }\omega,\\
\noame
\displaystyle \left(\int_0^{h(y')}\widetilde u'\,dy_3\right)\cdot n=0&\hbox{on }\partial\omega,\\
\noame
\displaystyle \widetilde u'=0&\hbox{on }\Gamma_0\cup \Gamma_1.
\end{array}\right.
\end{equation}

\end{theorem}
\begin{remark}
By uniqueness of solution of (\ref{limit_model}), we observe that the pair of functions $(\widetilde u, \widetilde p)$ are the same as those functions $(v^0,p^0)$ obtained in Section \ref{sec:formal} by formal arguments.
\end{remark}

\begin{proof}  We start by recalling that convergences (\ref{convergences_thm}) are consequence of Lemma \ref{lem_conv_vel} (it only remains to prove that $\partial_{y_3}\widetilde p_1=0$, which will be done at the end of Step 1). We divide the proof in two steps. In the first step, we derive the limit model  (\ref{limit_model}) and, in the second step, we prove the uniqueness of solution. The existence will be proved in the next corollary.\\

{\it Step 1. Derivation of the limit model.} Observe that conditions (\ref{limit_model})$_{2, 3, 4}$ are consequence of Lemma \ref{lem_conv_vel}. Now, let us pass to the limit in the variational inequality (\ref{v_ineq_carreau}) by taking into account (\ref{testv}).
\begin{itemize}
\item First term of (\ref{v_ineq_carreau}). From convergence (\ref{conv_vel_tilde}), since $\ep\mathbb{D}_{y'}[  \varphi']$ and $\varepsilon \mathbb{D}_{y'}[  { v}_\ep']$ tend to zero, we have
$$\begin{array}{l}
\displaystyle \eta_0\int_{\Omega}(1+\lambda|\ep\mathbb{D}_{y'}[  \varphi'] + \partial_{y_3}[ \varphi']|^2)^{{r\over 2}-1}(\varepsilon  \mathbb{D}_{y'}[  \varphi'] + \partial_{y_3}[ \varphi']):(\varepsilon \mathbb{D}_{y'}[  { v}_\ep'] + \partial_{y_3}[ { v}_\ep'])dy\\
\noame
=\displaystyle \eta_0\int_{\Omega}(1+\lambda|\partial_{y_3}[ \varphi']|^2)^{{r\over 2}-1}\partial_{y_3}[ \varphi']:\partial_{y_3}[\varphi'-\widetilde u']\,dy+O_\ep.
\end{array}$$
\item Second term of (\ref{v_ineq_carreau}). From the strong convergence of $p_\ep^0$ in $L^{r'}$ given in (\ref{pressuresr01})$_1$ and the weak convergence of velocity (\ref{conv_vel_tilde}), we have
$$\begin{array}{rl}\displaystyle 
\int_{\Omega}  p_\varepsilon^0\, {\rm div}_{y'}({v}'_\varepsilon)\,dy=&\displaystyle
\int_{\Omega}\widetilde p\, {\rm div}_{y'}(\varphi'-\widetilde u')\,dy+O_\varepsilon.
\end{array}$$

\item Third term of (\ref{v_ineq_carreau}). From ${\rm div}_{\ep}(\widetilde u_\ep)=0$ in $\Omega$, $\varphi_3\equiv 0$ and convergence (\ref{pressuresr01}) we have
$$\begin{array}{rl}\displaystyle 
\int_{\Omega}\widetilde p_\varepsilon^1\, {\rm div}_{\ep}({v}_\varepsilon)\,dy=&\displaystyle \int_{\Omega}\widetilde p_\varepsilon^1\, {\rm div}_{y'}(\varphi')\,dy\to 0.
\end{array}$$
\item Last term of (\ref{v_ineq_carreau}).  From convergence (\ref{conv_vel_tilde}), we have
$$\int_{\Omega} {f}'\cdot {v}'_\varepsilon\,dy= \int_{\Omega} {f}'\cdot (\varphi'-\widetilde{ u}')\,dy+O_\ep.$$
\end{itemize}
From previous convergences, we have the limit variational inequality
$$\begin{array}{l}\displaystyle
\eta_0\int_{\Omega}(1+\lambda|\partial_{y_3}[ \varphi']|^2)^{{r\over 2}-1}\partial_{y_3}[ \varphi']:\partial_{y_3}[ \varphi'-\widetilde { u}']\,dy
\displaystyle+\int_{\Omega}\widetilde p\, {\rm div}_{y'}(\varphi'-\widetilde u')\,dy\geq \int_{\Omega} {f}'\cdot (\varphi'-\widetilde{ u}')\,dy,
\end{array}$$
and by Minty's lemma is equivalent to
\begin{equation}\label{limit_var}\begin{array}{l}\displaystyle
{1\over 2}\eta_0\int_{\Omega}(1+{\lambda\over 2}|\partial_{y_3} \widetilde{ u}'|^2)^{{r\over 2}-1}\partial_{y_3}\widetilde { u}':\partial_{y_3}{  v}'\,dy
+\int_\Omega \widetilde p\,{\rm div}_{y'}(v')\,dy = \int_{\Omega} {f}'\cdot {v}'\,dy,
\end{array}
\end{equation}
which, by density, holds for every ${v}'\in V_{y_3}^r(\Omega)$. We observe that in (\ref{limit_var}), we have used that 
\begin{equation}\label{relationeq}\partial_{y_3}[\widetilde u']: \partial_{y_3} [{ v}']={1\over 2}\partial_{y_3}\widetilde{ u}' \cdot \partial_{y_3} { v}'\quad \hbox{and}\quad |\partial_{y_3}[\widetilde { u}']|^2={1\over 2}|\partial_{y_3}\widetilde{ u}'|^2.
\end{equation}
Integrating by parts, the limit variational formulation (\ref{limit_var}) is equivalent to the reduced Carreau Stokes system (\ref{limit_model}). The existence of solution of (\ref{limit_model}) will be proved in next corollary. Once the uniqueness of solution is proven (which will be done in the next step), then the whole sequence $(\ep^{-1}\widetilde u_\ep, \widetilde p_\ep)$ converges.

Finally, we comment that arguing as in the derivation of (\ref{v_ineq_carreau}) but with a test function of the type $\varphi=(\varphi',\varphi_3)$ with $\varphi'\equiv0$ and passing to the limit in this step, taking into account that $\widetilde u_3\equiv 0$, $\ep^{-1}\partial_{y_3}\widetilde p_1\rightharpoonup \partial_{y_3}\widetilde p_1$, and $f_3\equiv 0$, we can  derive $\partial_{y_3}\widetilde p_1=0$, which finishes the proof of (\ref{convergences_thm}).\\

 {\it Step 2. Uniqueness of solution. } We prove the uniqueness of solution $(\widetilde u',\widetilde p) \in V_{y_3}^r(\Omega)^3\times (L^{r'}_0(\omega)\cap W^{1,r'}(\omega))$ of (\ref{limit_model}) for $1<r<+\infty$ and $r\neq 2$.  The proof of uniqueness relies on the following inequalities satisfied for 
$$g_r(z')=\eta_0\left(1+{\lambda\over 2}|z'|^2\right)^{{r\over 2}-1}z',\quad \forall z'\in\mathbb{R}^2,$$
 given by (see for example \cite[Lemma 3.1]{Sandri} for more details)
 \begin{equation}\label{ineq_g}
 \begin{array}{l}
 (g_r(z')-g_r(t'),z'-t')_{\mathbb{R}^2}\geq   C_r |t'-z'|^2,\quad \hbox{if }r>2,
 \\
 \noame
 \displaystyle
 (g_r(z')-g_r(t'),z'-t')_{\mathbb{R}^2}\geq \overline C_r{|t'-z'|^2\over 1+|t'|^{2-r}+|z'|^{2-r}},\quad \hbox{if }1<r<2,
 \end{array}
 \end{equation}
 with $C_r,  \overline C_r>0$ independents of $t'$ and $z'$.\\

 Integrating (\ref{limit_model})$_1$ with respect to $y_3$, we obtain
 $$g_r(\partial_{y_3}\widetilde u')=C(y')-2y_3(f'(y')-\nabla_{y'}\widetilde p(y')).$$
 Let us suppose that (\ref{limit_model}) has two solutions $(\widetilde u_1', \widetilde p_1)$ and $(\widetilde u_2',\widetilde p_2)$, then 
 $$g(\partial_{y_3}\widetilde u_i')=C_i(y')-2y_3(f'(y')-\nabla_{y'}\widetilde p_i(y')),\quad i=1,2.$$
 For $1<r<2$ (the case $r>2$ is similar so we omit it), inequality (\ref{ineq_g})$_2$, says
 $$\left(g_r(\partial_{y_3}\widetilde u_2)-g_r(\partial_{y_3}\widetilde u_1),\partial_{y_3}(\widetilde u_2-\widetilde u_1)\right)\geq  \overline C_r{|\partial_{y_3}(\widetilde u_2-\widetilde u_1)|^2\over 1+|\partial_{y_3}\widetilde u_1|^{2-r}+|\partial_{y_3}\widetilde u_2|^{2-r}}.$$
Then, with $\alpha(y')=C_2(y')-C_1(y')$, we have
$$\left(\alpha(y')+2y_3\nabla_{y'}(\widetilde p_2-\widetilde p_1),\partial_{y_3}(\widetilde u_2-\widetilde u_1)\right)\geq   \overline C_r{|\partial_{y_3}(\widetilde u_2-\widetilde u_1)|^2\over 1+|\partial_{y_3}\widetilde u_1|^{2-r}+|\partial_{y_3}\widetilde u_2|^{2-r}}.$$
Integrating with respect to $y_3$, we get
$$\begin{array}{l}
\displaystyle \alpha(y')\int_0^{h(y')}\partial_{y_3}(\widetilde u_2-\widetilde u_1)\,dy_3+2\nabla_{y'}(\widetilde p_2-\widetilde p_1)\int_{0}^{h(y')}y_3\partial_{y_3}(\widetilde u_2-\widetilde u_1)\,dy_3\\
\noame
\displaystyle
\geq   \overline C_r\int_{0}^{h(y')}{|\partial_{y_3}(\widetilde u_2-\widetilde u_1)|^2\over 1+|\partial_{y_3}\widetilde u_1|^{2-r}+|\partial_{y_3}\widetilde u_2|^{2-r}}\,dy_3.
\end{array}$$
Due to the boundary conditions of $\widetilde u'$ given in (\ref{limit_model})$_4$, the first integral is zero, and, integrating by parts,
$$\begin{array}{l}
\displaystyle -2\nabla_{y'}(\widetilde p_2-\widetilde p_1)\int_{0}^{h(y')}\partial_{y_3}(\widetilde u_2-\widetilde u_1)\,dy_3+2\nabla_{y'}(\widetilde p_2-\widetilde p_1)[y_3(\widetilde u_2-\widetilde u_1]_{y_3=0}^{y_3=h(y')}\\
\noame
\displaystyle
\geq   \overline C_r\int_{0}^{h(y')}{|\partial_{y_3}(\widetilde u_2-\widetilde u_1)|^2\over 1+|\partial_{y_3}\widetilde u_1|^{2-r}+|\partial_{y_3}\widetilde u_2|^{2-r}}\,dy_3,
\end{array}$$
and again from (\ref{limit_model})$_4$, we get
$$\begin{array}{l}
\displaystyle -2\nabla_{y'}(\widetilde p_2-\widetilde p_1)\int_{0}^{h(y')}\partial_{y_3}(\widetilde u_2-\widetilde u_1)\,dy_3
\geq   \overline C_r\int_{0}^{h(y')}{|\partial_{y_3}(\widetilde u_2-\widetilde u_1)|^2\over 1+|\partial_{y_3}\widetilde u_1|^{2-r}+|\partial_{y_3}\widetilde u_2|^{2-r}}\,dy_3,
\end{array}$$
Integrating over $\omega$ and using Green's formula in the previous inequality, we get
$$\begin{array}{l}
\displaystyle\int_\omega (\widetilde p_2-\widetilde p_1)\left({\rm div}_{y'}\int_0^{h(y')}(\widetilde u_2-\widetilde u_1)\,dy_3\right)dy'-\int_{\partial\omega }(\widetilde p_2-\widetilde p_1)\left(n\cdot \int_0^{h(y')}(\widetilde u_2-\widetilde u_1)\,dy_3\right)dy' \,d\sigma\\
\noame
\displaystyle
\geq {1\over 2} \overline C_r\int_{\Omega}{|\partial_{y_3}(\widetilde u_2-\widetilde u_1)|^2\over 1+|\partial_{y_3}\widetilde u_1|^{2-r}+|\partial_{y_3}\widetilde u_2|^{2-r}}\,dy,
\end{array}$$
and with (\ref{limit_model})$_{2,3}$, we get
$$0\geq {1\over 2} \overline C_r\int_{\Omega}{|\partial_{y_3}(\widetilde u_2-\widetilde u_1)|^2\over 1+|\partial_{y_3}\widetilde u_1|^{2-r}+|\partial_{y_3}\widetilde u_2|^{2-r}}\,dy.$$
This implies $\partial_{y_3}\widetilde u_1=\partial_{y_3}\widetilde u_2$ a.e. in $\Omega$ and then, by using the boundary conditions (\ref{limit_model})$_4$, $\widetilde u_1=\widetilde u_2$. Moreover, 
$$\alpha(y')+2y_3\nabla_{y'}(\widetilde p_2(y')-\widetilde p_1(y'))=0,$$
which implies $\widetilde p_1=\widetilde p_2$ a.e. in $\omega$.

\end{proof}

In Corollary \ref{cor_thm}, we will need to consider the solution of an algebraic equation given in the following result, to give an expression for the Reynolds equation.
\begin{lemma}
The algebraic equation
\begin{equation}\label{taupsi}
\tau=\zeta\sqrt{{2\over \lambda}\left({\zeta\over \eta_0}\right)^{2\over r-2}-1},
\end{equation}
has a unique solution noted $\zeta=\psi(\tau)$ for $\tau\in \mathbb{R}^+$.
\end{lemma}
 \begin{proof}
 We follow \cite[Proposition 3.3]{Tapiero2} for the proof. Equation (\ref{taupsi}) can be written as follows
$$\zeta^2\left(\left({\zeta\over \eta_0}\right)^{2\over r-2}-1\right)={\lambda\over 2}\tau^2.$$
Function $\Theta:\zeta\to \zeta^2\left(\xi^{2\over r-2}-1\right)$ where $\xi =\zeta/\eta_0\geq 1$ is one to one from $[\eta_0,+\infty)$ into $[0,+\infty)$. Indeed, it holds
$$\Theta'(\zeta)=2\zeta\left[\xi^{2\over r-2}\left(1+{\zeta/\eta_0\over \xi(r-2)}\right)-1\right],$$
is strictly positive for $r>2$ and strictly negative for $1<r<2$, because
$$1+{\zeta/\eta_0\over \xi(r-2)}=1+{\xi\over \xi(r-2)}=1+{1\over r-2}={r-1\over r-2}<0.$$
\end{proof}

\begin{corollary}[Reynolds problem]\label{cor_thm}The Reynolds Velocity  
$$\widetilde V(y')=\int_0^{h(y')}\widetilde u(y)\,dy,$$
  is given by
\begin{equation}\label{Filtrationgamma1}\widetilde V'(y')=2\left( f'(y')-\nabla_{y'}\widetilde p(y')\right)\int_{-{h(y')\over 2}}^{h(y')\over 2}{({h(y')\over 2}+\xi)\xi \over \psi(2|f'(y')-\nabla_{y'}\widetilde p(y')||\xi|)}\,d\xi,\quad \widetilde V_3(y')=0,\quad\hbox{in }\omega.
\end{equation}
where $\psi$ is the inverse function of (\ref{taupsi}). The limit law for pressure can be written as follows
\begin{equation}\label{Reynoldsgamma1}
\left\{\begin{array}{l}
\displaystyle {\rm div}_{y'}\left(\left(\int_{-{h(y')\over 2}}^{h(y')\over 2}{({h(y')\over 2}+\xi)\xi \over \psi(2|f'(y')-\nabla_{y'}\widetilde p(y')||\xi|)}\,d\xi\right)\left(f'(y')-\nabla_{y'}\widetilde p(y')\right)\right)=0\quad \hbox{in }\omega,
\\
\noame
\displaystyle\left(\left(\int_{-{h(y')\over 2}}^{h(y')\over 2}{({h(y')\over 2}+\xi)\xi \over \psi(2|f'(y')-\nabla_{y'}\widetilde p(y')||\xi|)}\,d\xi\right)\left(f'(y')-\nabla_{y'}\widetilde p(y')\right)\right)\cdot n=0\quad \hbox{on }\partial\omega.
\end{array}\right.
\end{equation}

\end{corollary}
\begin{proof}  We follow   \cite[Proposition 3.3]{Tapiero2} for the proof. Writing (\ref{limit_model})$_1$ as follows
\begin{equation}\label{eqequiv}-\partial_{y_3}(\eta_r(\partial_{y_3}\widetilde u')\partial_{y_3}\widetilde u')=2(f'(y')-\nabla_{y'}\widetilde p(y')),\quad \hbox{with}\quad \eta_r(\xi')=\eta_0\left(1+{\lambda\over 2}|\xi'|^2\right)^{{r\over 2}-1},\quad \xi\in \mathbb{R}^2.
\end{equation}
we set $g(y')=2(f'(y')-\nabla_{y'}\widetilde p(y'))$ and integrating (\ref{eqequiv}) with respect to $y_3$, we deduce
\begin{equation}\label{eqequiv2}
\eta_r(\partial_{y_3}\widetilde u)\partial_{y_3}\widetilde u=C(y')-y_3 g(y').
\end{equation}
Set $\Gamma=C(y')-y_3 g(y')$ and $z=|\partial_{y_3}\widetilde u'|$ and so,  by (\ref{eqequiv}) we have that
$$\left(1+{\lambda\over 2}z^2\right)^{{r\over 2}-1}={\eta_r(z)\over \eta_0},$$
that is
\begin{equation}\label{eqequiv3}z=\sqrt{{2\over \lambda}\left\{\left({\eta_r(z)\over \eta_0}\right)^{2\over r-2}
-1\right\}}.
\end{equation}
Putting (\ref{eqequiv3}) in (\ref{eqequiv2}), we obtain 
$$\eta_r(z)\sqrt{{2\over \lambda}\left\{\left({\eta_r(z)\over \eta_0}\right)^{2\over r-2}
-1\right\}}=|\Gamma|,$$
which is (\ref{taupsi}) with $\zeta=\eta_r(z)$ and $\tau=|\Gamma|$.  Then,  we have $\eta_r(z)=\psi(|\Gamma|)$ and with (\ref{eqequiv2}), we deduce
$$\partial_{y_3}\widetilde u'={\Gamma\over \psi(|\Gamma|)}.$$ Integrating with respect to the vertical variable $y_3$ from $0$ to $y_3$ we get
\begin{equation}\label{uexpresion}
\widetilde u'=\int_0^{y_3}{\Gamma(y',\zeta)\over \psi(|\Gamma(y',\zeta)|)}\,d\zeta=\int_0^{y_3}{C(y')-\zeta g(y')\over \psi(|C(y')-\zeta g(y')|)}\,d\zeta.
\end{equation}
From the uniqueness of $\widetilde u'$, then constant $C(y')$ is uniquely determined. Indeed, if we use the boundary conditions of $\widetilde u'$ on $\Gamma_0\cup\Gamma_1$ given in (\ref{limit_model})$_4$, we have
\begin{equation}\label{eqequiv4}\int_0^{h(y')}{C(y')-\zeta g(y')\over \psi(|C(y')-\zeta g(y')|)}\,d\zeta=0.
\end{equation}
It holds that $C(y')={h(y')g(y')\over 2}$ because if we use this value of $C(y')$ and the change of variables $\xi={h(y')\over 2}-\zeta$ in (\ref{eqequiv4}), we have
$$\int_0^{h(y')}{g(y')({h(y')\over 2}-\zeta)\over \psi(|g(y')||{h(y')\over 2}-\zeta|)}\,d\zeta=g(y')\int_{-{h(y')\over 2}}^{h(y')\over 2}{\xi\over \psi(|g(y')||\xi|)}\,d\xi=0.$$
Then, we get expression
$$\widetilde u'(y)=g(y')\int_{0}^{y_3}{{h(y')\over 2}-\zeta\over \psi(|g(y')||{h(y')\over 2}-\zeta|)}\,d\zeta=g(y')\int_{{h(y')\over 2}-y_3}^{h(y')\over 2}{\xi\over \psi(|g(y')||\xi|)}\,d\xi,$$
which proves the existence of  solution $\widetilde u'$, which is unique, of (\ref{limit_model}).

Finally, integrating $\widetilde u'$ with respect to $y_3$ between $0$ and $h(y')$, we get
$$\begin{array}{rl}\displaystyle 
\widetilde V'(y')=\int_0^{h(y')}\widetilde u'(y)\,dy_3= &\displaystyle g(y')\int_0^{h(y')}\left(\int_{{h(y')\over 2}-y_3}^{h(y')\over 2}{\xi\over \psi(|g(y')||\xi|)}\,d\xi\right)\,dy_3\,\\
\noame
=&\displaystyle g(y')\int_{-{h(y')\over 2}}^{h(y')\over 2}{\xi\over \psi(|g(y')||\xi|)}\left(\int_{{h(y')\over 2}-\xi}^{h(y')}\,dy_3\right)\,d\xi,
\end{array}$$
that is, expression 
$$\widetilde V'(y')=\displaystyle g(y')\int_{-{h(y')\over 2}}^{h(y')\over 2}{({h(y')\over 2}+\xi)\xi \over \psi(|g(y')||\xi|)}d\xi.
$$
 We remark that, due to $\widetilde u_3=0$, we have $\widetilde V_3=0$.

Finally, by using expression (\ref{Filtrationgamma1}) of $\widetilde V'$ and divergence conditions (\ref{limit_model})$_{2,3}$, we deduce the Reynolds problem (\ref{Reynoldsgamma1}) for $\widetilde p$.
 \\

\end{proof}

 \section*{Conflict of interest} The authors state that there is no conflict of interest to declare.

 \section*{Data availability statement}
Data sharing not applicable to this article as no datasets were generated or analysed during the current study.

\end{document}